\documentclass[12pt,letterpaper]{article}

\pdfoutput=1

\usepackage{physics}
\usepackage{epsfig}
\usepackage{amsmath}
\usepackage{amssymb}
\usepackage{amsthm}
\usepackage{indentfirst}
\usepackage{xspace}
\usepackage{multirow}
\usepackage{hyperref}
\usepackage{xcolor}
\usepackage{verbatim}
\usepackage{subfigure}
\usepackage{mathrsfs}
\usepackage{ytableau}
\usepackage{bbm}

\usepackage[letterpaper,margin=0.9in,headheight=15pt]{geometry}
\usepackage{mathpazo}
\usepackage{authblk}
\usepackage{empheq}
\usepackage{feynmp}
\usepackage{graphicx}
\usepackage[matrix,arrow]{xy}
\usepackage{young}
\usepackage[vcentermath]{youngtab}
\usepackage{slashed}
\usepackage{bbm}
\usepackage{youngtab}
\usepackage{rotfloat}
\usepackage{stmaryrd}
\usepackage{amsfonts,amssymb,amsmath}
\usepackage{tikz-cd}
\usepackage{thmtools}
\usepackage{dashrule}
\usepackage[missing=]{gitinfo2}
\usepackage{fancyhdr}
\usepackage{mdframed}
\usepackage{booktabs}
\usepackage{subfiles}
\usepackage{simplewick}

\usepackage[utf8]{inputenc}

\definecolor{darkblue}{rgb}{0.1,0.1,0.7}
\definecolor{darkred}{rgb}{0.5,0.1,0.1}
\definecolor{darkgreen}{rgb}{0.0,0.42,0.06}
\hypersetup{colorlinks=true,urlcolor=darkblue,linkcolor=darkblue,citecolor=darkred}
\definecolor{shadecolor}{rgb}{0.85,0.85,0.85}

\declaretheoremstyle[spaceabove=0.25cm,spacebelow=0.25cm,notefont=\normalfont\bfseries, notebraces={(}{)}]{theorem}
\declaretheoremstyle[spaceabove=0.25cm,spacebelow=0.25cm,bodyfont=\normalfont,notefont=\normalfont\bfseries, notebraces={(}{)}]{noital}
\declaretheoremstyle[spaceabove=0.25cm,spacebelow=0.25cm,bodyfont=\normalfont\color{darkgreen},notefont=\normalfont\bfseries, notebraces={(}{)}]{green}
\declaretheoremstyle[spaceabove=0.25cm,spacebelow=0.25cm,bodyfont=\normalfont,notefont=\normalfont\bfseries,qed=$\qedsymbol$,notebraces={(}{)}]{proofstyle}

\numberwithin{equation}{section}

\usepackage[most,listings]{tcolorbox}
\tcbuselibrary{breakable}

\tcbset{colframe=white}

\newcommand{\bC}{\ensuremath{\mathbb{C}}}

\newcommand{\bZ}{\ensuremath{\mathbb{Z}}}

\newcommand{\Tw}{\mathbbm{T}\!\mathbbm{w}}
\newcommand{\BB}{\mathbbm{B}}
\newcommand{\DD}{\mathbbm{D}}
\newcommand{\FF}{\mathbbm{F}}

\def \be  {\begin{equation}}
\def \ee  {\end{equation}}
\def \bea {\begin{equation}\begin{aligned}}
\def \eea {\end{aligned}\end{equation}}
\def \ba  {\begin{eqnarray}}
\def \ea  {\end{eqnarray}}

\newtheorem{lemma}{Lemma}[section]
\newtheorem{conjecture}[lemma]{Conjecture}

\def\ybox{\Yboxdim{5pt}\yng(1)}

\begin{document}

\title{Cyclotomic expansions of HOMFLY-PT  colored by rectangular Young diagrams}
\author{Masaya Kameyama, Satoshi Nawata, Runkai Tao, Hao Derrick Zhang}
\date{}

\maketitle

{\abstract{We conjecture a closed-form expression of HOMFLY-PT invariants of double twist knots colored by rectangular Young diagrams where the twist is encoded in interpolation Macdonald polynomials. We also put forth a conjecture of cyclotomic expansions of HOMFLY-PT polynomials colored by rectangular Young diagrams for any knot.}}


\tableofcontents

\section{Introduction}

Colored HOMFLY-PT polynomials are two-variable quantum knot invariants associated with irreducible representations of  $A_N$ type Lie algebras.
Although several methods to compute HOMFLY-PT polynomials for arbitrary color in principle are known, carrying out explicit computation is practically very
challenging for general non-torus knots and colors. However, in recent years, studying the structural properties of colored HOMFLY-PT polynomials, closed-form expressions for symmetric representations have been found for a certain class of non-torus knots.

Actually, the structural properties become more apparent at the level of HOMFLY-PT homology that categorifies quantum HOMFLY-PT polynomials.
Lately, the HOMFLY-PT homology colored by arbitrary representations has been defined in \cite{cautis2016remarks}. Although it is formidable to carry out computation of homology via the definition, various structural properties of HOMFLY-PT homology have been uncovered by combining mathematical definitions and physical predictions. In particular, it was proposed in \cite{Gorsky:2013jxa,Gukov:2015gmm} that when colors are specified by rectangular Young diagrams, structural properties become more manifest if $(\mathbf{a,Q,t_r,t_c})$ quadruple-grading is introduced.

\begin{figure}
  \begin{minipage}[b]{8cm}\centering
	\centering
	\raisebox{-3cm}{\includegraphics[width=6cm]{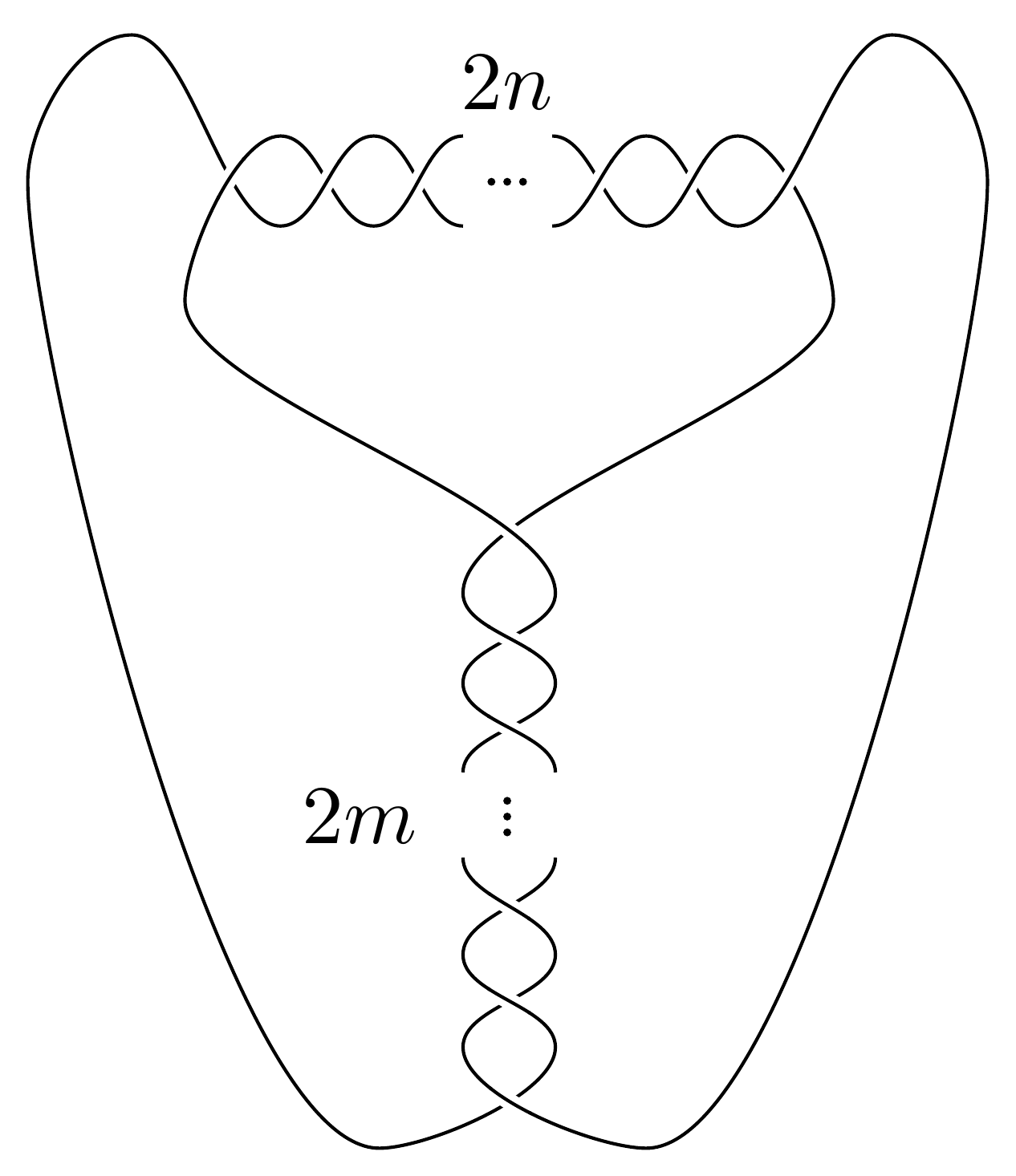}}
	\label{fig:two-braid-knot}
	\end{minipage}
	  \begin{minipage}[b]{8cm}\centering
	 $$ \begin{array}{c | c | c}m&n& \textrm{Name}\\ \hline{1}&{1}&{3 _{1}} \\{- 1}&{1}&{4 _{1}} \\{2}&{1}&{5 _{2}} \\{- 2}&{1}&{6 _{1}} \\{3}&{1}&{7 _{2}} \\{2}&{2}&{7 _{4}} \\{- 3}&{1}&{8 _{1}} \\{- 2}&{2}&{8 _{3}} \\{4}&{1}&{9 _{2}}\\{3}&{2}&{9 _{5}}  \\{- 4}&{1}&{10 _{1}} \\{- 3}&{2}&{10 _{3}} \end{array}$$
	  \end{minipage}
	 \caption{Double twist knot $K_{m,n}$ and corresponding knots  in Rolfsen Table up to 10 crossings.}\label{fig:double-twist-knot}
\end{figure}

In \cite{Kononov:2016cwp}, closed-form expressions of Poincare polynomials of HOMFLY-PT homology of the figure-eight and the trefoil colored by rectangular Young diagrams $[r^s]$ are conjectured. The conjectural formulas in   \cite{Kononov:2016cwp} are very simple, being expressed by a summation over Young diagrams inscribed by $[r^s]$. In this paper, we further generalize the formulas to the case of double twist knots drawn in Figure \ref{fig:double-twist-knot}. Our conjectural  formulas are expressed by a summation  over iteratively inscribed Young diagrams with interpolation Macdonald polynomials, which can be understood as a generalization of cyclotomic expansions  from symmetric representations \cite{habiro2008unified,Nawata:2012pg,Gukov:2015gmm}  to rectangular Young diagrams. We attach Mathematica files to arXiv page, which explicitly compute both HOMFLY-PT polynomials and Poincar\'e polynomials of the double twist knots colored by rectangular Young diagrams. Like the Rosso-Jones formula \cite{Rosso:1993vn} for torus knots, we believe that the formulas will find many applications in other areas of mathematics.

\subsubsection*{Convention}
Throughout this paper, we use the following skein relation for a \emph{reduced} HOMFLY polynomial
${H}_{\ybox}(K;A,q)$:
\bea\nonumber
&&A ~{H}_{\ybox}\left({\raisebox{-.2cm}{\includegraphics[width=.6cm]{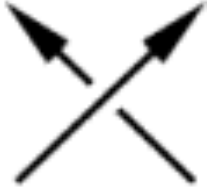}}}\right)
- A^{-1}{H}_{\ybox}\left({\raisebox{-.2cm}{\includegraphics[width=.6cm]{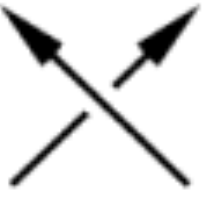}}}\right)
=
(q-q^{-1}) {H}_{\ybox}\left({\raisebox{-.2cm}{\includegraphics[width=.6cm]{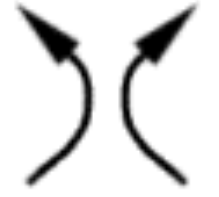}}}\right)\,,
\eea
with the unknot invariant is normalized as
\be\nonumber
{H}_{\ybox} ({\raisebox{-.1cm}{\includegraphics[width=.4cm]{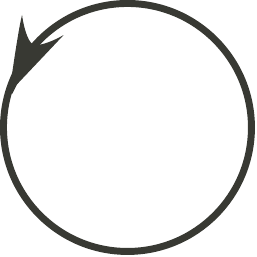}}}) = 1~. \ee
In addition, in this paper, a knot is always zero-framed and we do not consider non-trivial framings.

\subsubsection*{Acknowledgement}
We would like to thank Alexei Morozov for discussion and correspondence. Without his help and encouragement, this paper would be impossible. S.N. would like to thank ICTS, KITP and  IHES for warm hospitality where a part of work has been carried out. R.T. and H.D.Z. are grateful for the support from NSFC grant No.11850410428.  This research was supported  in part by the International Centre for Theoretical Sciences (ICTS) during a visit for participating in the program - Quantum Fields, Geometry and Representation Theory (Code: ICTS/qftgrt2018/07), and in part by the National Science Foundation under Grant No. NSF PHY-1748958 during a visit  for participating in  Quantum Knot Invariants and Supersymmetric Gauge Theories.

\section{Colored HOMFLY-PT polynomials}

It is  conjectured in \cite{Kononov:2016cwp} that the $[r^s]$-colored HOMFLY-PT polynomial of the figure-eight $K_{1,-1}$ is expressed as
$$
H_{[r^s]}(K_{1,-1};A,q) =
\sum_{\lambda\subset [r^s]}
D_\lambda^{(s)}(q) D_{\lambda^T}^{(r)}(q^{-1})\prod_{\ybox \in \lambda}
\{
Aq^{r + a'_{\ybox}-l'_{\ybox}}
\}
\{
A q^{a'_{\ybox} -l'_{\ybox}-s}
\}\, ,
$$
where $\{x\} =  x-x^{-1}$ and $\lambda^T$ is the transposition of $\lambda$. Here, the factor $D_\lambda^{(s)}(q)$ is the principal specialization of the Schur polynomial
$$
D_\lambda^{(s)}(q):=s_\lambda(q^{s-1},q^{s-3},\ldots,q^{-s+3},q^{-s+1} )=\prod_{\ybox \in \lambda}\frac{\{
	q^{s-l'_{\ybox}+a'_{\ybox}}
	\}    }
{	\{q^{l_{\ybox} + a_{\ybox}+1}\}
}~.
$$
The formula is deduced by repackaging tedious calculations \cite{Morozov:2016eqp} with the Cauchy formula
$$
\prod _{i , j}\left( 1 + x _{i}y _{i}\right) = \sum _{\lambda}s _{\lambda}( x ) s _{\lambda^T }( y )~,
$$
so that it naturally exhibits the exponential growth property of colored HOMFLY-PT polynomials \cite{Zhu:2012tm} with respect to colors
$$
H_\lambda(K;A,q=1)=\Big[H_{\ybox}(K;A,q=1)\Big]^{|\lambda|}~,
$$
where $|\lambda|$ is the number of boxes in $\lambda$.
It is an elegant extension from symmetric representations \cite{Itoyama:2012fq} to rectangular Young diagrams. Hence, as in symmetric representations \cite{Nawata:2012pg,Gukov:2015gmm}, one can expect that  the $[r^s]$-colored HOMFLY-PT polynomial of the double twist knot can be expressed by inserting an appropriate twist element into this formula.

\begin{figure}
  \begin{minipage}[b]{8cm}\centering
	\centering
\includegraphics[width=5cm]{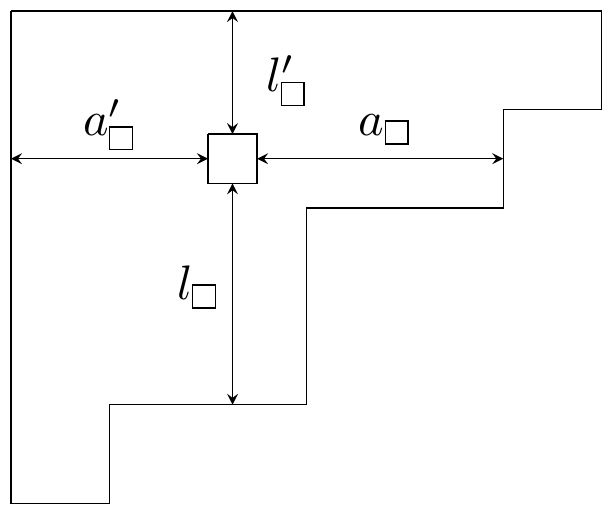}
	 \caption{Convention}
	\label{fig:Young-convention}
	\end{minipage}
	  \begin{minipage}[b]{8cm}\centering
\includegraphics[width=6cm]{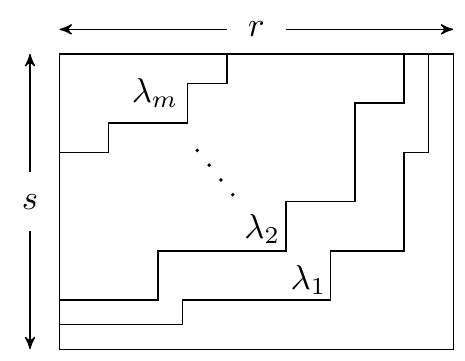}
	 \caption{Iteratively inscribed Young diagrams $[r^s]\supset\lambda_1\supset\cdots\supset \lambda_m$.}\label{fig:iterative}
	  \end{minipage}
\end{figure}

In the case of symmetric representations, the twist element has been found from the $q$-binomial theorem.
In fact, the binomial formula for any representation $\lambda$ is given in  \cite[\S3, Example 10]{macdonald1998symmetric}:
\be\label{Schur-1+x}
s _{\lambda}\left( 1 + x _{1}, 1 + x _{2}, \ldots , 1 + x _{n}\right) = \sum _{\mu \subset \lambda} d _{\lambda \mu}s _{\mu}\left( x _{1}, \ldots , x _{n}\right)
\ee
where
$$
d _{\lambda \mu}= \operatorname{det}\left( \left( \begin{array}{c}{ \lambda _{i}+ n - i}\\{\mu _{j}+ n - j}\end{array} \right) \right) _{1 < i , j < n}~.
$$
Thus, we use its $q$-deformation to define
\begin{equation}\label{B-factor}
B_{\lambda,\mu}^{(s)}(A,q) := A^{2\abs{\mu}}
q^{(1-s)(\abs{\lambda}-\abs{\mu})}
\qty(\prod_{\ybox \in \mu}
q^{4(a'_{\ybox} - l'_{\ybox})})\  \ \textrm{det}
\qty(
{
\lambda_i + s -i \brack
\mu_j + s -j
}_{q^2}
)_{1\leq i, j \leq s}\, .
\end{equation}
where ${x \brack y}_q$ is a $q$-binomial coefficient. Note that when $\mu=\emptyset$, we have $B_{\lambda,\emptyset}^{(s)}=D_\lambda^{(s)}$.
Using this, we introduce the following factor which admits two expressions
\bea\label{F-unrefined}
F_{\lambda}^{(m)}(A,q)&=\qty(D_\lambda^{(s)})^{-1}\sum_{\lambda=\lambda_1\supset\cdots\supset\lambda_m}
 B_{\lambda_1,\lambda_2}^{(s)} B_{\lambda_2,\lambda_3}^{(s)}
\cdots B_{\lambda_{m-1},\lambda_{m}}^{(s)}D_{\lambda_{m}}^{(s)}
\cr
&=\qty(D_{\lambda^T}^{(r)})^{-1}\sum_{{\lambda^T}=\lambda'_1\supset\cdots\supset\lambda'_m}
 B_{\lambda'_1,\lambda'_2}^{(r)} B_{\lambda'_2,\lambda'_3}^{(r)}
\cdots B_{\lambda'_{m-1},\lambda'_{m}}^{(r)}D_{\lambda'_{m}}^{(r)}\Bigg|_{q\to q^{-1}}
\,,
\eea
where the summations are taken over iteratively inscribed Young diagrams as illustrated in Figure \ref{fig:iterative}.   Although the building blocks $B^{(s)}$ and $D^{(s)}$  depend on $s$ in the first line, $F_{\lambda}^{(m)}(A,q)$ is a Laurent polynomial of $(A,q)$  \emph{independent of} $s (\ge \textrm{length}(\lambda))$. In a similar fashion, $F_{\lambda}^{(m)}(A,q)$ is independent of  $r(\ge \textrm{length}(\lambda^T))$ though each building block depends on $r$ in the second line.
Note that $F^{m=1}_\lambda=1$ for any $\lambda$.

Using this, we define the twist element as
\begin{align}
  \textrm{Tw}^{(m)}_\lambda(A,q): = \begin{cases}(-A^2)^{\abs{\lambda}}\Big(\displaystyle\prod_{\ybox \in \lambda}q^{2 a_{\ybox} - 2 l_{\ybox}} \Big) F^{(m)}_\lambda(A,q) \, , &  m\in \bZ_{>0}   \, ,\\
    F^{(-m)}_\lambda( A^{-1},  q^{-1})\, , & m\in \bZ_{<0} \, .\end{cases}
\end{align}
Then, we conjecture that  the HOMFLY-PT polynomial of the double twist knot $K_{m,n}$ colored by a rectangular Young diagram $[r^s]$ is expressed as
\begin{equation}
H_{[r^s]}(K_{m,n};A,q) =
\sum_{\lambda\subset [r^s]} \frac{\textrm{Tw}^{(m)}_\lambda \textrm{Tw}^{(n)}_\lambda}
{\textrm{Tw}^{(1)}_\lambda \textrm{Tw}^{(-1)}_\lambda}
D_\lambda^{(s)}(q) D_{\lambda^T}^{(r)}(q^{-1})\prod_{\ybox \in \lambda}
\{
Aq^{r + a'_{\ybox}-l'_{\ybox}}
\}
\{
A q^{a'_{\ybox} -l'_{\ybox}-s}
\}\, .
\end{equation}
This is a natural extensions of the formulas for symmetric representations \cite{Nawata:2012pg,Gukov:2015gmm}  to rectangular Young diagrams.
We have checked this formula with the results in \cite{Kononov:2016hle,Morozov:2017ezc,Morozov:2017lvo}.

\bigskip

Motivated by this formula, we further generalize the conjecture of cyclotomic expansions of colored HOMFLY-PT polynomials of a knot $K$ for symmetric representations \cite[Eqn.(2)]{Kononov:2015dda}, \cite[Conjecture 2.3]{Nawata:2015wya} to rectangular Young diagrams.
\begin{conjecture}[Cyclotomic expansions]\label{HOMFLY-cyclotomic}
Let $Y$ be a set of all Young diagrams. For a knot $K$, there exists a function
 \be
C(K): Y \times \bZ_{>0}\longrightarrow\bZ[A^\pm,q^\pm] \ ; \ (\lambda,s)\mapsto C_{\lambda,s}(K;A,q)~,\nonumber
\ee
which satisfies the following properties
\begin{itemize}\setlength{\parskip}{-0.1cm}
\item $C_{\emptyset,s}(K;A,q)=1$ for any $s\in \bZ_{>0}$
\item $C_{\lambda,s}(K;A,q)=0$ for $s<\textrm{length}(\lambda)$
\end{itemize}
such that the $[r^s]$-colored reduced HOMFLY-PT polynomial of the knot $K$ is expressed as
\bea\label{cyclotomic-HOMFLY}
H_{[r^s]}(K;A,q)
&=&(A^{\bullet}q^{\bullet})^{rs} \sum_{\lambda\subset [r^s]} C_{\lambda,s}(K;A,q)  D_{\lambda^T}^{(r)}(q^{-1})\prod_{\ybox \in \lambda}
\{
Aq^{r + a'_{\ybox}-l'_{\ybox}}
\}
~.
\eea
Importantly, $C_{\lambda,s}(A,q)$ is independent of $r$.
\end{conjecture}
In fact, when a color is a symmetric reperesentation $\lambda=[r]$, \eqref{cyclotomic-HOMFLY} reduces to
$$
H_{[r]}(K ; A, q)=\left(A^{\bullet} q^{\bullet}\right)^{r} \sum_{r \geq k \geq 0} C_{k}(K ; A, q) q^{-2 r k}(A^{2} q^{2 r} ; q^{2})_{k}
{r  \brack k}_{q^{2}}~,
$$
which is exactly equal to \cite[Conjecture 2.3]{Nawata:2015wya}.

If we can make use of the transposition symmetry  \cite{Zhu:2012tm} of the colored HOMFLY-PT polynomials
$$
H_\lambda(K;A,q)=H_{\lambda^T}(K;A,q^{-1})~,
$$
then the conjecture can be stated by exchanging $q\to q^{-1}$, $s\to r$ and $\lambda\to \lambda^T$, but the details are omitted.

\section{Poincar\'e polynomials of colored HOMFLY-PT homology}

In \cite{Kononov:2016cwp},  it is further conjectured that Poincar\'e polynomial of quadruply-graded HOMFLY-PT homology of the figure-eight colored by $[r^s]$ is given by
$$
\mathscr{P}_{[r^s]}(K_{1,-1};\mathbf{a,Q,t_r,t_c}) = \sum_{\lambda\subset [r^s]}
\DD_\lambda^{(s)}(q,t) \DD_{\lambda^T}^{(r)}(t^{-1},q^{-1})\prod_{\ybox \in \lambda}
\{ \frac A \sigma q^{r + a'_{_\square}} t^{- l'_{\square}}\} \{{A}{\sigma}  q^{  a'_{_\square}} t^{-s  - l'_{_\square}}\} \bigg|_{\eqref{cov}}\,
$$
where $\DD_\lambda^{(s)}$ is a principal specialization of the Macdonald polynomial
$$
\DD_\lambda^{(s)}(q,t):=P_\lambda(t^{s-1},t^{s-3},\ldots,t^{-s+3},t^{-s+1} ;q,t)=\prod_{\ybox \in \lambda}\frac{\{t^{s-l'_{_\square}} q^{a'_{_\square}}\}}{\{t^{l_{_\square}+1} q^{a_{_\square}}\}}~,
$$
and the change of variables is
\be\label{cov}A\to \mathbf{a \sqrt{-t_r t_c}}, \quad t\to \mathbf{t_r}^{-1} ,\quad q\to - \mathbf{t_c}, \quad \sigma\to \mathbf{t_r}^{-s}\mathbf{Q}^{-1}
~.\ee

To find the corresponding twist element, we need the binomial theorem involving Macdonald polynomials. Happily, the binomial theorem has been generalized by using interpolation Macdonald polynomials $P_\lambda^*(x,q,t)$ \cite{okounkov1997binomial}. Some basics facts on this matter are summarized in Appendix \ref{sec:app-a}.
Hence, motivated by the $(q,t)$-version \eqref{qt-binomial} of the binomial theorem, we define the refined version of \eqref{B-factor}
\bea
\BB_{\lambda,\mu}^{(s)}(A,q,t): =&(-t^{(1-s)})^{\abs{\lambda}-\abs{\mu}} \Big(\frac{Aq}{t}\Big)^{2\abs{\mu}}\Big(\displaystyle\prod_{\ybox \in \lambda}q^{2 a'_{\ybox} }t^{-2 l'_{\ybox}} \Big) \Big(\displaystyle\prod_{\ybox \in \mu}q^{2a'_{\ybox} }t^{- 2l'_{\ybox}} \Big) \cr& \times \frac{P _ {\mu}^ {*}( q^{2\lambda} ; q^{2},t^{2})}{P _ {\mu}^ {*}( q^{2\mu} ; q^{2},t^{2})}.
\frac{P _ {\lambda}^ {*}( 0 ; q^{-2},t^{-2})}{P _ {\mu}^ {*}(0; q^{-2},t^{-2})}~.
\eea
In general, $\BB_{\lambda,\mu}^{(s)}$ is not a Laurent polynomial but a rational function, and $\BB_{\lambda,\emptyset}^{(s)}=\DD_\lambda^{(s)}$.
Then, we introduce the factor
\bea
\FF_{\lambda}^{(m)}(A,q,t)&=\qty(\DD_\lambda^{(s)})^{-1}\sum_{\lambda=\lambda_1\supset\cdots\supset\lambda_m}
 \BB_{\lambda_1,\lambda_2}^{(s)} \BB_{\lambda_2,\lambda_3}^{(s)}
\cdots \BB_{\lambda_{m-1},\lambda_{m}}^{(s)}\DD_{\lambda_{m}}^{(s)}\cr
&=\qty(\DD_{\lambda^T}^{(r)})^{-1}\sum_{{\lambda^T}=\lambda'_1\supset\cdots\supset\lambda'_m}
 \BB_{\lambda'_1,\lambda'_2}^{(r)} \BB_{\lambda'_2,\lambda'_3}^{(r)}
\cdots \BB_{\lambda'_{m-1},\lambda'_{m}}^{(r)}\DD_{\lambda'_{m}}^{(r)}\Bigg|_{q\leftrightarrow t^{-1}}
\,.
\eea
As in \eqref{F-unrefined}, it is independent of $s~  (\ge\textrm{length}(\lambda))$ and $r~ (\ge \textrm{length}(\lambda^T))$. Remarkably, even though each building block is a rational function of $(A,q,t)$, $\FF_{\lambda}^{(m)}$ is always a Laurent polynomial of  $(A,q,t)$. We also observe that if $\textrm{length}(\lambda)\le2$, it is a Laurent polynomial with positive coefficients whereas it generally involves both positive and negative coefficients for $\textrm{length}(\lambda)>2$.

Finally, we conjecture that the Poincar\'e polynomial of $[r^s]$-colored HOMFLY-PT homology of the double twist knot $K_{m,n}$ can be expressed as
\bea\label{final}
&\mathscr{P}_{[r^s]}(K_{m,n};\mathbf{a,Q,t_r,t_c}) =(-\mathbf{t_r}\mathbf{t_c})^{-rs(\textrm{sgn}(m)+\textrm{sgn}(n))/2}\times \cr
&\sum_{\lambda\subset [r^s]}
\DD_\lambda^{(s)}(q,t) \DD_{\lambda^T}^{(r)}(t^{-1},q^{-1})\prod_{\ybox \in \lambda}
\{ \frac A \sigma q^{r + a'_{_\square}} t^{- l'_{\square}}\} \{{A}{\sigma}  q^{  a'_{_\square}} t^{-s  - l'_{_\square}}\} \bigg|_{\eqref{cov}}\,  \frac{\Tw^{(m)}_\lambda \Tw^{(n)}_\lambda}
{\Tw^{(1)}_\lambda \Tw^{(-1)}_\lambda}~,
\eea
where the twist elements are defined as
$$
\Tw^{(m)}_\lambda(\mathbf{a},\mathbf{t_r},\mathbf{t_c}) := \begin{cases}\Big(-\mathbf{a}^2\mathbf{t_r}^2\mathbf{t_c}^2\Big)^{\abs{\lambda}}
\Big(\displaystyle\prod_{\ybox \in \lambda}\mathbf{t_c}^{2 a_{\ybox} }\mathbf{t_r}^{ 2 l_{\ybox}} \Big) \FF^{(m)}_\lambda(\mathbf{a},  \mathbf{t_c},\mathbf{t_r}) \, , &  m\in \bZ_{>0}   \, ,\\
    \FF^{(-m)}_\lambda(\mathbf{a}^{-1},  \mathbf{t_c}^{-1},\mathbf{t_r}^{-1})\, , & m\in \bZ_{<0} \, .\end{cases}
$$
Surprisingly, $\mathscr{P}_{[r^s]}(K_{m,n})$ is a Laurent polynomial of $(\mathbf{a,Q,t_r,t_c})$ with \emph{positive coefficients} even though summands in \eqref{final} are rational functions of $(\mathbf{a,Q,t_r,t_c})$ in general. We have checked the structural properties \cite{Gorsky:2013jxa} such as self-symmetry, mirror symmetry, refined exponential growth property and colored differentials
for a number of examples.

\section{Discussion}
The formulas conjectured in this paper cry out for geometric interpretation. In fact, the binomial theorem \eqref{Schur-1+x} can be interpreted by the first Chern class of the tensor product of vector bundles over a flag variety \cite{lascoux1982classes}. Recently, a new definition of the uncolored HOMFLY-PT homology is given by $(q,t)$-equivariant sheaves on Hilbert schemes of points on $\bC^2$  \cite{gorsky2016flag,oblomkov2018knot,oblomkov2018categorical}.  Although its colored version has yet to be defined, the formula in this paper strongly suggests that braiding on $(q,t)$-equivariant sheaves will be captured by interpolation Macdonald polynomials via the localization of equivariant Grothendieck-Riemann-Roch formula \cite{haiman2002notes} if it is defined.

The double twist knots can be obtained by taking surgeries on two components of Borromean rings with framings $-1/m$ and $-1/n$.
Actually, the cyclotomic expansions of colored Jones polynomials of the double twist knots have been originally obtained from that of Borromean rings \cite{habiro2008unified}. Therefore the structure in this formula can be extended to Borromean rings and twist links at the level of colored HOMFLY-PT polynomials as in \cite{habiro2008unified,Gukov:2015gmm}.
In addition, it would be interesting to extract information about quantum $6j$-symbols with rectangular Young diagrams from the formula conjectured in this paper.

Although thirty years have passed since colored quantum knot invariants have been introduced, our understanding of colored invariants beyond symmetric representations is still very limited.
The formula in this paper reveals deep mathematical structure hidden behind colored knot invariants, relating to special functions and potentially geometry. Of course, we just glimpse a tip of the iceberg, and we hope
that our results will serve as a stepping stone toward the study of knot invariants of general colors.

\appendix
\section{Interpolation Macdonald polynomials}
\label{sec:app-a}

In this ``Appendix'', we review a combinatorial definition of interpolation Macdonald polynomials introduced by \cite{sahi1996interpolation,knop1997symmetric,okounkov1998shifted} and their binomial theorem \cite{okounkov1997binomial}.

For a Young diagram $\lambda$, a Young tableau $T$ is obtained by filling the boxes of the Young diagram by the numbers in $\{1, ..., n\}$.
We denote the entry of a box by $T_{\ybox}$.
A Young tableau is called a \emph{reverse tableau} if the entries is weakly decreasing along each row and strongly decreasing along each column.
The maximal number of entries in a reverse tableau must be not less than the length of the Young diagram by definition. A reverse tableau $T$ on $\lambda$ makes a sequence of Young diagrams
$$
 \emptyset\equiv \lambda^{(n)}\subseteq \lambda^{(n-1)} \subseteq \cdots \subseteq \lambda^{(0)}\equiv \lambda,
$$
where $\lambda^{(k)}$ is the shape of a sub Young tableau of $T$ such that all entries satisfy $T_{\ybox}>k$.
For example, if we consider the following reverse tableau on $\lambda=(4,4,3,2,2)$ with $n=6$
\begin{center}
{
\Yboxdim15pt
\young(6662,5521,332,22,11)~,
}
\end{center}
then the sequence is
$$
{\Yboxdim6pt
\emptyset \subseteq \yng(3) \subseteq \yng(3,2) \subseteq \yng(3,2) \subseteq \yng(3,2,2) \subseteq \yng(4,3,3,2) \subseteq \yng(4,4,3,2,2)~.
}
$$
A skew Young diagram $\lambda/\mu$  is called \emph{horizontal strip} if there is at most one box in each column of $\lambda$.
For a horizontal strip $\lambda/\mu$, we denote sets of rows and columns intersecting with $\lambda/\mu$ by $R_{\lambda/\mu}$ and $C_{\lambda/\mu}$. Then, $R_{\lambda/\mu}-C_{\lambda/\mu}$ is a set of boxes belonging to $R_{\lambda/\mu}$ but not to $C_{\lambda/\mu}$. For the above sequence $\lambda^{(k-1)}/\lambda^{(k)}$ are horizontal strips for every $k=1,\ldots,6$. As an example, $\lambda^{(0)}/\lambda^{(1)}$ and $R_{\lambda^{(0)}/\lambda^{(1)}}-C_{\lambda^{(0)}/\lambda^{(1)}}$ correspond to the black boxes and the red box, respectively
\begin{center}
\ytableausetup{boxsize=13pt}
  \begin{ytableau}
   *(white) &   & & \\
     &   &*(red)  &*(black) \\
    & & \\
&  \\
*(black) & *(black)
  \end{ytableau}
\end{center}

Let $x=(x_1,\ldots,x_n)$ be $n$-tuple variables.
The ordinary Macdonald polynomial $P_{\lambda}(x;q,t)$ is combinatorially defined as
$$
 P_{\lambda}(x;q,t) = \sum_{T} \psi_{T}(q,t) \prod _ {\ybox \in \lambda}  x_{T_{\ybox}}~,
$$
where the sum is taken over all reverse tableaux $T$ on $\lambda$ with entries in $\{1,\ldots,n\}$ and the coefficients $\psi_{T}(q,t)$ are defined by
$$
 \psi _{T}( q,t)  = \prod _{i = 1}^{n}\psi _{\lambda ^{( i - 1 )}/ \lambda ^{( i )}} ( q,t) ~,
\quad
\psi _{\lambda / \mu}( q,t) = \prod _{\ybox \in R_{ \lambda / \mu}- C_{ \lambda / \mu}} \frac{b _{\mu}( \ybox )}{ b _{\lambda}( \ybox  )}~,
$$
with
$$
 b _{\lambda}( \ybox )  = \frac{1 - q ^{a _{\ybox}}t ^{l _{\ybox}+ 1}}{1 - q ^{a _{\ybox}+ 1}t ^{l _ {\ybox}}}~.
$$
Then, the interpolation Macdonald polynomial $P _ {\mu}^{*}( x ; q,t)$ is combinatorially defined by
$$
P _ {\lambda}^{*}( x ; q,t) = \sum _ {T}\psi _ {T}( q,t) \prod _ {\ybox \in \mu}t ^{1 - T _{\ybox}}\left( x _ {T _{\ybox}}- q ^{a ^{\prime}_{\ybox}}t ^{- l ^{\prime}_{\ybox}}\right).
$$
Using the interpolation Macdonald polynomials, the binomial theorem is generalized \cite{okounkov1997binomial} to
$$
\frac{P _{\lambda}^{*}\left(cx _{1}, \ldots ,cx _{n}; q,t\right)}{ P _{\lambda}^{*}(c, \ldots ,c; q,t)}= \sum _{\mu}\frac {c^{| \mu |}}{t ^{( n - 1 ) | \mu |}} \frac{P _{\mu}^{*}\left( q ^{- \lambda}; q^{-1} , t^{-1} \right)}{ P _{\mu}^{*}\left( q ^{- \mu}; q^{-1} , t^{-1} \right)}\frac{P _{\mu}^{*}\left( x _{n}, \ldots , x _{1}; q^{-1} , t^{-1} \right)}{ P _{\mu}^{*}(c, \ldots ,c; q,t) }~.
$$
In particular, the limit $c\to 0$ leads to
\be\label{qt-binomial}
\frac{P _{\lambda}^{*}\left( x _{1}, \ldots , x _{n}; q,t\right)}{ P _{\lambda}^{*}( 0 , \ldots , 0 ; q,t)}= \sum _{\mu}\frac{P _{\mu}^{*}\left( q ^{- \lambda}; q^{-1} , t^{-1} \right)}{ P _{\mu}^{*}\left( q ^{- \mu}; q^{-1} , t^{-1} \right)}\frac{P _{\mu}\left( x _{1}, x _{2} t^{-1}, \ldots , x _{n}t ^{1 - n}; q,t\right)}{ P _{\mu}^{*}( 0 , \ldots , 0 ; q,t) }~.
\ee
If we restrict ourselves to the cases of symmetric representations and single variable $x_1$, this formula reduces to the usual $q$-binomial theorem.

\newcommand{\etalchar}[1]{$^{#1}$}
\providecommand{\bysame}{\leavevmode\hbox to3em{\hrulefill}\thinspace}
\providecommand{\MR}{\relax\ifhmode\unskip\space\fi MR }
\providecommand{\MRhref}[2]{%
  \href{http://www.ams.org/mathscinet-getitem?mr=#1}{#2}
}
\providecommand{\href}[2]{#2}

\bigskip

\bigskip

\noindent
Masaya Kameyama, \textsf{Graduate School of Mathematics, Nagoya University,  Nagoya 464-8602, Japan},
{\href{mailto:m13020v@math.nagoya-u.ac.jp}{m13020v@math.nagoya-u.ac.jp}

\

\noindent
Satoshi Nawata, \textsf{Department of Physics and Center for Field Theory and Particle Physics, Fudan University, 2005
Songhu Road, 200438  Shanghai, China}\\
\textsf{Institut des Hautes \'Etudes Scientifiques,
35 Route de Chartres, Bures-sur-Yvette, 91440, France} \\
\textsf{Kavli Institute for Theoretical Physics, Santa Barbara CA 93106, U.S.A.}, \href{mailto:snawata@gmail.com}{snawata@gmail.com}

\

\noindent
Runkai Tao, \textsf{Department of Physics and Center for Field Theory and Particle Physics, Fudan University, 2005
Songhu Road, 200438 Shanghai, China}, \href{mailto:runkaitao@gmail.com}{runkaitao@gmail.com}

\

\noindent
Hao Derrick Zhang, \textsf{Department of Physics and Center for Field Theory and Particle Physics, Fudan University, 2005
Songhu Road, 200438 Shanghai, China}, \href{mailto:haozhangphys@gmail.com}{haozhangphys@gmail.com}

\end{document}